\documentclass[12pt,psfig,reqno]{amsart}
\usepackage{amsfonts}
\usepackage{amscd}
\usepackage{amssymb,amsfonts,latexsym}
\usepackage{graphics,verbatim}
\usepackage{graphicx}
\usepackage[usenames]{color} %\color{Red}
\makeatletter
\renewcommand\section{\@startsection {section}{1}{\z@}%
                                   {-3.5ex \@plus -1ex \@minus -.2ex}%
                                   {2.3ex \@plus.2ex}%
                                   {\centering\normalfont\bf}}

 \numberwithin{equation}{section}
\numberwithin{equation}{section}

\numberwithin{equation}{section}

\theoremstyle{plain}

\begin{document}
\title{ Composition Operators on $\bf H^{p,q,s}(B_{n})$ of $\bf \mathbb{C}^{n}$ }
\author{Hongxin Chen$^{1,*}$ , \ \  Xuejun Zhang$^{2}$ }
\address{1. School of Mathematics,  Hunan  University, Changsha, Hunan, China \\
2. College of Mathematics and Statistics, Hunan Normal University, Changsha, Hunan, China   }

\email{1755310775@qq.com; \ xuejunttt@263.net}

\date{}
\keywords {Composition operator;   boundedness; compactness;  general Hardy type space;  unit ball} \subjclass[2010]{
32A37; 47B33
 }
\thanks{$^*$ Corresponding author.\\
 The research is supported by the Natural Science Foundation of Hunan Province of
 China (No. 2022JJ30369) and the Education
Department Important Foundation of Hunan Province in China (No. 23A0095)}

\begin{abstract}   Let $B_{n}$ be the unit ball in the complex vector space $\mathbb{C}^{n}$, and  let $\varphi: B_{n}\rightarrow B_{n}$ be a holomorphic mapping. In this paper, we characterize those symbols $\varphi$ such that composition operators $C_{\varphi}$ are bounded or compact  on the general Hardy type space
$H^{p,q,s}(B_{n})$. These results extend the relevant results on Hardy space and  some other classical function spaces.
\end{abstract}
\maketitle
\section{\bf Introduction\label{sect.1}}

\ \ \ Suppose that $\mathbb{C}^{n}$ is the $n$-dimensional complex vector space. For $w=(w_{1},\cdots,w_{n})$ and $u=(u_{1},\cdots,u_{n})$ in $\mathbb{C}^{n}$, let
 $\langle w, u\rangle= w_{1}\overline{u_{1}}+\cdots +w_{n}\overline{u_{n}}$,  and let $B_{n}=\{w=(w_{1},\cdots, w_{n})\in \mathbb{C}^{n}: |w|=\sqrt{\langle w,w\rangle}<1\}$ be the  unit ball in $\mathbb{C}^{n}$. Suppose that   $H(B_{n})$ is the set of all holomorphic functions on $B_{n}$, and  $H^{\infty}(B_{n})$ is the set  of
all bounded holomorphic functions on $B_{n}$.
Let $dv$ be the normalized volume measure on $B_{n}$,  and $d\sigma$ be the normalized rotation invariant  measure
on  $ S_{n}$, where $S_{n}$ is the boundary of $B_{n}$.

\vskip2mm

 {\bf Definition 1.1} \ For $p>0$,  if \ $h\in H(B_{n})$ and
 $$||h||_{p}=\sup_{u\in B_{n}}(1-|u|^{2})^{p}\ |
h(u)|<\infty, $$
then  $h$ is said to belong to  the holomorphic growth  space
 $H^{\infty}_{p}(B_{n})$.
\vskip2mm
{\bf Definition 1.2} \ For $p>0$ and $q> 0$, the weighted Hardy space $H_{q}^{p}(B_{n})$ consists of all holomorphic functions   $h$  on $B_{n}$ such that
$$
||h||_{p,q}=\sup_{0\leq r<1}\left\{(1-r^{2})^{q}\int_{S_{n}}|h(r\eta)|^{p}\
d\sigma(\eta)\right\}^{\frac{1}{p}}<\infty.
$$

{\bf Definition 1.3} \ For $\beta>-1$ and $p>0$, the function $h$ is said to belong to
 the weighted Bergman space $A_{\beta}^{p}(B)$ if $h\in H(B_{n})$ and
$$
||h||_{A_{\beta}^{p}}=\left(\int_{B}|h(u)|^{p}\
dv_{\beta}(u)\right)^{\frac{1}{p}}<\infty,
$$
where $dv_{\beta}(u)=c_{\beta}(1-|u|^{2})^{\beta}\ dv(u)$, the
constant $c_{\beta}$ such that $v_{\beta}(B_{n})=1$.
\vskip2mm

For  $w\in B_{n}$, let $\varphi_{w}$ be the holomorphic automorphism of $B_{n}$  with $\varphi_{w}(0) = w$,  $\varphi_{w}(w)=0$
and $\varphi_{w}^{-1}=\varphi_{w}$.
\vskip2mm

  {\bf Definition 1.4} \ For  $p>0$,  $s\geq 0$, $q+n\geq
0$,  $q+s\geq 0$,  if $f$ is  a Lebesgue  measurable function on $B_{n}$
and $||f||_{p,q,s}=\displaystyle{\sup_{0\leq r<1}}M_{p,q,s}(r,f)<\infty$,
 where
\begin{align*}
 M_{p,q,s}(r,f)=\sup_{a\in B_{n}}\left\{(1-r^{2})^{q}\int_{S_{n}}|f(r\xi)|^{p}\
(1-|\varphi_{a}(r\xi)|^{2})^{s}\ d\sigma(\xi)\right\}^{\frac{1}{p}},
\end{align*}
then $f$ is said to belong to the general Lebesgue space $L^{p,q,s}(B_{n})$.
\vskip2mm
 The space $L^{p,q,s}(B_{n})$ is a Banach space under the norm
$||. ||_{p,q,s}$ when $p\geq 1$. If $0<p<1$, then  $L^{p,q,s}(B_{n})$ is
a complete metric space under the  distance
$\rho(f,g)=||f-g||^{p}_{p,q,s}$.
In particular,  $H^{p,q,s}(B_{n})=L^{p,q,s}(B_{n})\bigcap H(B_{n})$ is called
the general Hardy type space. In fact, $H^{p,q,s}(B_{n})$  is just the Hardy space
$H^{p}(B_{n})$ when $q=s=0$.   Therefore,   $H^{p,q,s}(B_{n})$  is a generalization  of $H^{p}(B_{n})$.
The space $H^{p,q,s}(B_{n})$ comes from some practical applications (for example, [1-3]).
 In [4], we formally defined this function class as the general Hardy type space on the unit ball.  The space
$H^{p,q,s}(B_{n})$ contains several
   classical function spaces by taking different parameters $q$ and $s$, for examples,  $H^{\infty}(B_{n})$ ($s\geq n=-q$),
    $H^{\infty}_{\frac{q+n}{p}}(B_{n})$ ($s\geq n>-q$),  $H^{p}_{q}(B_{n})$ ($q>0=s$) etc. If $0<s<n$, then  $H^{p,q,s}(B_{n})$
    is a  new holomorphic function space,  different from the Hardy space $H^{k}(B_{n})$ and the weight Hardy space $H^{k}_{t}(B_{n})$ for any $0<k\leq \infty$ and $0<t<\infty$ (see [4]).

\vskip2mm

{\bf Definition 1.5} \ Let  $\varphi: B_{n}\rightarrow B_{n}$ be a holomorphic mapping. The $\varphi$  induces a operator $C_{\varphi}$ on $H(B_{n})$: $C_{\varphi}g=g\circ\varphi$ \ \ ($g\in H(B_{n}$)).  The operator $C_{\varphi}$ is called  composition operator.
\vskip1mm

 Composition type  operators  between various  holomorphic function spaces
 on  the unit disc  or ball have been
studied for a long time, and there have been a lot of research results (for examples, [1]-[23]).
The main purpose of this paper is to discuss the boundedness and compactness of $C_{\varphi}$ on $H^{p,q,s}(B_{n})$. These results extend the relevant results on Hardy space and  some other classical function spaces.
\vskip2mm
If there exists constant $c>0$  such that $ E_{1} \geq cE_{2}$
 ( or $ E_{1} \leq cE_{2}$), then we write $``E_{1} \gtrsim E_{2}" $  ( or $``E_{1}\lesssim E_{2}"$). If $``E_{1} \gtrsim E_{2}" $ and $``E_{1}\lesssim E_{2}"$, then we call $``E_{1}\asymp E_{2}"$.

\section{\bf  Some  Lemmas \label{sect.2}}

\ \ \ {\bf Lemma 2.1 ([4])} \ For $p>0$, $s\geq 0$, $q+s\geq 0$, $q+n\geq 0$, if \ $h\in H^{p,q,s}(B_{n})$, then
$$
|h(z)|\lesssim\frac{||h||_{p,q,s}}{(1-|z|^{2})^{\frac{q+n}{p}}} \ \
\mbox{for all $z\in B_{n}$}.
$$

This result comes from Lemma 2.3 in [4].

\vskip2mm
{\bf Lemma 2.2 ([24])} \ For $\beta>-1$, if $h\in A_{\beta}^{1}(B_{n})$, then
$$
h(z)=\int_{B_{n}}\frac{h(w) \ dv_{\beta}(w)}{(1-\langle z, w\rangle)^{n+1+\beta}} \ \ (z\in B_{n}).
$$

This result comes from Theorem 2.2 in [24].

\vskip2mm

{ \bf Lemma 2.3 ([26])} \ Let  $c>0$ and $\delta>-1$. Then for all $w\in B_{n}$,
$$
\int_{S_{n}}\frac{d\sigma(\xi)}{|1-\langle
w,\xi\rangle|^{n+c}}\asymp \int_{
B_{n}}\frac{(1-|u|^{2})^{\delta}\ dv(u)}{|1-\langle
w,u\rangle|^{n+1+\delta+c}}\asymp\frac{1}{(1-|w|^{2})^{c}}.
$$

This result comes from  Proposition 1.4.10 in [26].

\vskip2mm For two points in $B_{n}$, there are the following
integral bidirectional estimates. These results come from Proposition 3.1 in [25].
 \vskip2mm {\bf Lemma 2.4 ([25])} \  For  $t> 0$ and $\delta> 0$, \ let
$$
\displaystyle{I_{w,u}=\int_{S_{n}}\frac{d\sigma(\eta)}{|1-\langle
\eta,w\rangle|^{ \delta}\ |1-\langle \eta,u\rangle|^{ l}} \ \ (w,u\in B_{n}).}$$

Then there are the following bidirectional estimates:

\vskip2mm

(1) \ \ $ \displaystyle{I_{w,u}\asymp \frac{1}{  |1-\langle w,u\rangle|^{\delta+l-n}}\ }$  when $\delta+l>n$ and  $\max\{\delta,l\}<n$.
\vskip2mm

(2) \ \ $ \displaystyle{I_{w,u}\asymp\frac{1}{|1-\langle w,u\rangle|^{l}}\log\frac{e}{|1-\langle w,\varphi_{w}(u)\rangle|} }
$ when $\delta=n>l$.
\vskip2mm

(3) \ \ $ \displaystyle{I_{w,u}\asymp \frac{1}{|1-\langle w,u\rangle|^{n}}\log\frac{e}{1-|\varphi_{w}(u)|^{2}}} $  when $\delta=n=l$.
\vskip2mm

(4) \ \ $\displaystyle{I_{w,u} \asymp\frac{(1-|w|^{2})^{
n-\delta}}{|1-\langle w,u\rangle|^{ l}} + \frac{(1-|u|^{2})^{ n-l}}{
|1-\langle w,u\rangle|^{ \delta}}}$  \ when $\delta>n$ and $ l>n$.

\vskip2mm

(5) \ \ $\displaystyle{I_{w,u} \asymp \frac{(1-|w|^{2})^{
n-\delta}}{|1-\langle w,u\rangle|^{ n}} + \frac{\log\frac{e}{1-|\varphi_{u}(w)|^{2}}}{|1-\langle
w,u\rangle|^{ \delta}}}$ \ when
$\delta>n=l$.

\vskip2mm

 (6) \ \
$ \displaystyle{I_{w,u}\asymp \frac{1}{(1-|w|^{2})^{ \delta-n}\
|1-\langle w,u\rangle|^{ l}}} $ \ when \  $\delta>n>l$.

\vskip3mm {\bf Lemma 2.5 ([23])} \ For $p>0$ and $0<p_{0}\leq 1$, \ then
 $$
 (1-\rho^{2})^{\frac{(1-p_{0})n}{p_{0}}}\int_{S_{n}}|h(\rho^{2}\zeta)|^{p}\ d\sigma(\zeta)\leq
\left\{\int_{S_{n}}|h(\rho\zeta)|^{pp_{0}}\
d\sigma(\zeta)\right\}^{\frac{1}{p_{0}}}
 $$
 for all $h\in H(B_{n})$ and $0\leq \rho<1$.

\vskip2mm
The result comes from Lemma 2.3 in [23].

\section{\bf   Main  Results \label{sect.3}}

\ \ \ If $s\geq n$, then $H^{p,q,s}(B_{n})=B^{\frac{q+n}{p}}(B_{n})$, the  Bloch type space. If $s=0$, then $H^{p,q,s}(B_{n})$ is the weighted Hardy space or the Hardy space. The boundedness and compactness of $C_{\varphi}$ on the Bloch type space or the weighted Hardy space  have been solved.  Therefore, we only need to consider  the case $0< s<n$.
 \vskip2mm
 {\bf Theorem 3.1} \ Let \ $p>0$, \ $0<s<n$, \ $q+s\geq 0$, \ $q+n\geq 0$.  Suppose that \ $\varphi: \ B_{n}\rightarrow B_{n}$ is a holomorphic mapping.
 \vskip2mm

 (1) \ If \ $\varphi$ is an automorphism, then  $C_{\varphi}$ is  bounded  on $H^{p,q,s}(B_{n})$.
 \vskip2mm

 (2) \ If \ $q+n>0$ and $C_{\varphi}$ is  compact  on $H^{p,q,s}(B_{n})$, then
  \begin{align}
 \lim_{|w|\rightarrow 1^{-}}\frac{1-|w|^{2}}{1-|\varphi(w)|^{2}}=0.
 \end{align}

(3) \ If \ $||\varphi||_{\infty}<1$, then  $C_{\varphi}$ is  compact on $H^{p,q,s}(B_{n})$.

\vskip2mm
(4) \ If $C_{\varphi}$ is  bounded  on $H^{p,q,s}(B_{n})$, then for any $t>0$ there is
\begin{align}
\sup_{0\leq \rho<1}\sup_{ a, w\in B_{n}}(1-\rho^{2})^{q}\int_{S_{n}}\frac{(1-|w|^{2})^{t}(1-|\varphi_{a}(\rho \xi)|^{2})^{s}}{|1-\langle \varphi(\rho\xi),w\rangle|^{q+n+t}}\ d\sigma(\xi)<\infty.
\end{align}

(5) \ When $q+s>0$, for any  $\beta > \{q+\max(s,pn)\}\min(1/p,1)-1$, \ if
\begin{align}
\sup_{\eta\in S_{n}}\int_{S_{n}}\frac{(1-|a|^{2})^{s}|1-\langle\varphi(\rho \xi),\varphi(a)\rangle|^{2s}\ d\sigma(\xi)}{(1-|\varphi(a)|^{2})^{s}|1-\langle \varphi(\rho \xi), r\eta\rangle|^{n+1+\beta}|1-\langle\rho \xi,a\rangle|^{2s}}\lesssim \frac{1}{(1-r\rho)^{\beta+1}}
\end{align}
for all $a\in B_{n}$ and $0\leq r,\rho<1$, then $C_{\varphi}$ is  bounded  on $H^{p,q,s}(B_{n})$.
\vskip2mm
(6) \ When $q+s>0$, if (3.3) holds, then $C_{\varphi}$ is a compact operator on $H^{p,q,s}(B_{n})$  if and only if (3.1) holds.
\vskip2mm

 {\bf Proof} \  For the following proof, we select constant  $\alpha>\max\{(q+n)/p-1, \ (q+s)/p-1, \ (q+s+n)/p-n-1\}$.
 \vskip2mm

 (1) \ This result comes from Theorem 5.1 in [31].

\vskip2mm

(2) \ If $||\varphi||_{\infty}<1$, then  (3.1) is obvious. When $||\varphi||_{\infty}=1$, let $\{w^{j}\}\subset B_{n}$ be  any sequence such that
$|\varphi(w^{j})|\rightarrow 1$ when $j\rightarrow\infty$. We take
$$
g_{j}(w)=\frac{(1-|\varphi(w^{j})|^{2})^{\frac{|q|}{p}+1}}{(1- \langle w,\varphi(w^{j})\rangle)^{\frac{q+n+|q|}{p}+1}} \ \ \ (w\in B_{n}).
$$
Then $\{g_{j}\}$ converges to 0 uniformly on any compact set of $B_{n}$. Otherwise,
it follows from Lemma 2.4(4-6) and $\displaystyle{\sup_{0<x\leq1}x^{\varepsilon}\log\frac{e}{x}}<\infty$ ($\varepsilon>0$) (case $2s=n$) that
$$
\sup_{ 0\leq \rho<1}\sup_{a\in B_{n}}\int_{S_{n}}\frac{(1-r^{2})^{q+s}(1-|a|^{2})^{s}(1-|\varphi(w^{j})|^{2})^{|q|+p}}{|1-\langle r\xi, \varphi(w^{j})\rangle|^{q+n+|q|+p}|1-\langle r\xi,a\rangle|^{2s}}\ d\sigma(\xi)\lesssim 1
$$
for all  $j=1,2,\cdots$ by means of three cases $0\leq 2s<n$, $2s=n$ and $2s>n$. Therefore,  it follows from Lemma 2.1 that
\begin{align}
&
\left(\frac{1-|w^{j }|^{2}}{1-|\varphi(w^{j})|^{2}}\right)^{\frac{q+n}{p}}=(1-|w^{j}|^{2})^{\frac{q+n}{p}}|C_{\varphi}g_{j}(w^{j})|\nonumber\\
&\lesssim ||C_{\varphi}g_{j}||_{p,q,s}\rightarrow 0 \ \ (j\rightarrow\infty).
\end{align}
This implies that (3.1) holds.

\vskip2mm
(3) \ When $||\varphi||_{\infty}<1$, for any sequence $\{h_{k}\}$ which converges to 0 uniformly on any compact set of $B_{n}$ and $\displaystyle{\sup_{k\in \{1,2,\cdots\}}}||h_{k}||_{p,q,s}\leq 1$, it follows from Lemma 2.3 that
\begin{eqnarray*}
 &\;&
\sup_{0\leq \rho<1}\sup_{a\in B_{n}}(1-\rho^{2})^{q}\int_{S_{n}}|h_{k}[\varphi(\rho \xi)]|^{p}(1-|\varphi_{a}(\rho \xi)|^{2})^{s}\ d\sigma(\xi)\\
&\;&
\leq \max_{|w|\leq ||\varphi||_{\infty}}|h_{k}(w)|^{p}\sup_{0\leq \rho<1}\sup_{a\in B_{n}}\int_{S_{n}}\frac{(1-\rho^{2})^{q+s}(1-|a|^{2})^{s}\ d\sigma(\xi)}{|1-\langle \rho\xi,a \rangle|^{2s}} \\
&\;&
\lesssim \max_{|w|\leq ||\varphi||_{\infty}}|h_{k}(w)|^{p} \ \rightarrow 0 \ \ (k\rightarrow\infty) \ \Rightarrow \ \mbox{$C_{\varphi}$ is compact on $H^{p,q,s}(B_{n})$.}
\end{eqnarray*}

(4) \ For any $w\in B_{n}$ and $t>0$, let
$$
f_{w}(z)=\frac{(1-|w|^{2})^{\frac{t}{p}}}{(1-\langle z, w\rangle)^{\frac{q+n+t}{p}}} \ \ (z\in B_{n}).
$$

 By Lemma 2.4(1-6) and $\displaystyle{\sup_{0<x< 2}x^{\varepsilon}\log\frac{e}{x}}<\infty$ ($\varepsilon>0$),  we have that
 $$
 \sup_{0\leq \rho<1}\sup_{a\in B_{n}}\int_{S_{n}}\frac{(1-\rho^{2})^{q+s}(1-|a|^{2})^{s}(1-|w|^{2})^{t}}{|1-\langle \rho\xi, w\rangle|^{q+n+t}|1-\langle \rho\xi,a\rangle|^{2s}}\ d\sigma(\xi)\lesssim 1
 $$
 for all $w\in B_{n}$. It follows from the boundedness of $C_{\varphi}$ on $H^{p,q,s}(B_{n})$
that
\begin{eqnarray*}
 &\;&
||C_{\varphi}||\gtrsim ||C_{\varphi}||.||f_{w}||_{p,q,s}\geq ||C_{\varphi}f_{w}||_{p,q,s}\\
&\;&
=\sup_{0\leq \rho<1}\sup_{a\in B_{n}}(1-\rho^{2})^{q}\int_{S_{n}}\frac{(1-|w|^{2})^{t}(1-|\varphi_{a}(\rho \xi)|^{2})^{s}}{|1-\langle \varphi(\rho\xi),w\rangle|^{q+n+t}}\ d\sigma(\xi)
\end{eqnarray*}
 for all $w\in B_{n}$. This shows that (3.2) holds.

\vskip2mm
(5) \ For any $g\in H^{p,q,s}(B_{n})$ and $a\in B_{n}$, we take
$$
G_{a,g}(w)=\frac{g(w)}{(1-\langle w,\varphi(a)\rangle)^{\frac{2s}{p}}} \ \ (w\in B_{n}).
$$
It follows from  $\alpha>(q+n)/p-1$ and Lemmas 2.1-2.2 that
$$
G_{a,g}(w)=\int_{B_{n}}\frac{G_{a,g}(u)\ dv_{\alpha}(u)}{(1-\langle w, u)^{n+1+\alpha}} \ \ (w\in B_{n}).
$$
This means that
\begin{align}
&
|G_{a,g}[\varphi(z)]|\leq \int_{B_{n}}\frac{|G_{a,g}(u)|\ dv_{\alpha}(u)}{|1-\langle \varphi(z), u\rangle|^{n+1+\alpha}}\nonumber\\
&=\int_{B_{n}}\frac{|1-\langle u,\varphi(a)\rangle|^{-\frac{2s}{p}}|g(u)|\ dv_{\alpha}(u)}{|1-\langle \varphi(z), u\rangle|^{n+1+\alpha}} \ \ (z\in B_{n}).
\end{align}

If \ $p>1$, \ then  we may choose $\max\{\alpha, \ q+s-1\}<\alpha_{1}<\min\{p(\alpha+1)-1, \ q+s+\alpha\}$ such that $(p\alpha-\alpha_{1})/(p-1)>-1$ and $\alpha_{1}-\alpha>0$.  By  (3.5), H\"{o}lder inequality,  Lemma 2.3, $(1-|z|^{2})/(1-|\varphi(z)|^{2})\leq (1+|\varphi(0)|)/(1-|\varphi(0)|)$, we have that
\begin{align}
&
|G_{a,g}[\varphi(z)]|^{p}
 \lesssim \int_{B_{n}}\frac{(1-|u|^{2})^{\alpha_{1}}|g(u)|^{p}\ dv(u)}{|1-\langle u,\varphi(a)\rangle|^{2s}|1-\langle \varphi(z), u\rangle|^{n+1+\alpha}}\nonumber\\
 &
 \times\left\{\int_{B_{n}}\frac{(1-|u|^{2})^{\frac{p\alpha-\alpha_{1}}{p-1}}\ dv(u)}{|1-\langle \varphi(z), u\rangle|^{n+1+\alpha}}\right\}^{p-1}\nonumber\\
 &
 \asymp (1-|\varphi(z)|^{2})^{\alpha-\alpha_{1}}\int_{B_{n}}\frac{(1-|u|^{2})^{\alpha_{1}}|g(u)|^{p}\ dv(u)}{|1-\langle u,\varphi(a)\rangle|^{2s}|1-\langle \varphi(z), u\rangle|^{n+1+\alpha}} \nonumber\\
&
\lesssim (1-|z|^{2})^{\alpha-\alpha_{1}}\int_{B_{n}}\frac{(1-|u|^{2})^{\alpha_{1}}|g(u)|^{p}\ dv(u)}{|1-\langle u,\varphi(a)\rangle|^{2s}|1-\langle \varphi(z), u\rangle|^{n+1+\alpha}}.
\end{align}

It follows from the selection of $\alpha$ and $\alpha_{1}$ that there are  $\alpha_{1}-q-s>-1$ and $q+s+\alpha-\alpha_{1}>0$.   When (3.3) holds, for any $0\leq \rho<1$ and $\xi\in S_{n}$, by (3.6),
    Lemma 1.8 in [24], (3.3), Lemma 6 in [27], we have that
\begin{eqnarray*}
 &\;&
(1-\rho^{2})^{q}\int_{S_{n}}|g[\varphi(\rho \xi)]|^{p}(1-|\varphi_{a}(\rho \xi)|^{2})^{s}\ d\sigma(\xi)\\
 &\;&
 =(1-\rho^{2})^{q+s}(1-|a|^{2})^{s}\int_{S_{n}}|G_{a,g}[\varphi(\rho \xi)]|^{p}\frac{|1-\langle \varphi(\rho\xi),\varphi(a)\rangle|^{2s}}{|1-\langle \rho\xi,a\rangle|^{2s}}\ d\sigma(\xi)
 \\
 &\;&
\lesssim \int_{S_{n}}\int_{B_{n}}\frac{(1-\rho^{2})^{q+s+\alpha-\alpha_{1}}|g(u)|^{p}|1-\langle\varphi(\rho\xi),\varphi(a)\rangle|^{2s}(1-|u|^{2})^{\alpha_{1}}\ dv(u)d\sigma(\xi)}{
(1-|a|^{2})^{-s}|1-\langle\varphi(\rho\xi), u\rangle|^{n+1+\alpha}|1-\langle\rho\xi,a\rangle|^{2s}|1-\langle u,\varphi(a)\rangle|^{2s}}\\
 &\;&
 =\int_{B_{n}}\frac{|g(u)|^{p}(1-|u|^{2})^{\alpha_{1}}}{|1-\langle u,\varphi(a)\rangle|^{2s}}\int_{S_{n}}\frac{(1-\rho^{2})^{q+s+\alpha-\alpha_{1}}
 |1-\langle\varphi(\rho\xi), \varphi(a)\rangle|^{2s}\ d\sigma(\xi) dv(u)}{(1-|a|^{2})^{-s}|1-\langle\varphi(\rho \xi), u\rangle|^{n+1+\alpha}|1-\langle\rho\xi,a\rangle|^{2s}}\\
 &\;&
\lesssim \int_{0}^{1}(1-r^{2})^{\alpha_{1}-q-s}\int_{S_{n}}(1-r^{2})^{q}|g(r\eta)|^{p}(1-|\varphi_{\varphi(a)}(r \eta)|^{2})^{s}\\
 &\;&
\times \left\{\frac{(1-|a|^{2})^{s}}{(1-|\varphi(a)|^{2})^{s}}\int_{S_{n}}\frac{(1-\rho^{2})^{q+s+\alpha-\alpha_{1}}|1-\langle\varphi(\rho\xi),\varphi(a)\rangle|^{2s}\ d\sigma(\xi)}{|1-\langle\varphi(\rho \xi),r\eta|^{n+1+\alpha}|1-\langle\rho\xi,a\rangle|^{2s}}\right\} d\sigma(\eta)dr\\
 &\;&
\lesssim (1-\rho)^{q+s+\alpha-\alpha_{1}}||g||_{p,q,s}^{p}\int_{0}^{1}\frac{(1-r)^{\alpha_{1}-q-s}}{(1-r\rho)^{\alpha+1}}\ dr
\asymp ||g||_{p,q,s}^{p}.
\end{eqnarray*}
This means that $C_{\varphi}$ is a bounded operator on $H^{p,q,s}(B_{n})$.
\vskip2mm
If \ $0<p\leq 1$, \ then we let $\alpha=(\alpha'+n+1)/p-n-1$. It follows from $\alpha>(q+s+n)/p-n-1$ and $\alpha>(q+n)/p-1$ that $\alpha'>q+\max\{s,pn\}-1$.
When (3.3) holds,  by (3.5), Lemma 2.15 in [24],
 Fubini theorem,  Lemma 1.8 in [24], (3.3),  Lemma 6 in [27], we have that
\begin{eqnarray*}
 &\;&
(1-\rho^{2})^{q}\int_{S_{n}}|g[\varphi(\rho\xi)]|^{p}(1-|\varphi_{a}(\rho \xi)|^{2})^{s}\ d\sigma(\xi)
 \\
 &\;&
\lesssim \int_{S_{n}}\int_{B_{n}}\frac{(1-\rho^{2})^{q+s}|g(u)|^{p}|1-\langle\varphi(\rho \xi),\varphi(a)\rangle|^{2s}(1-|u|^{2})^{\alpha'}\ dv(u)d\sigma(\xi)}{(1-|a|^{2})^{-s}|1-\langle\varphi(\rho \xi),u\rangle|^{n+1+\alpha'}|1-\langle\rho \xi,a\rangle|^{2s}|1-\langle u,\varphi(a)\rangle|^{2s}}\\
 &\;&
 =\int_{B_{n}}\frac{|g(u)|^{p}(1-|u|^{2})^{\alpha'}}{|1-\langle u,\varphi(a)\rangle|^{2s}}\int_{S_{n}}\frac{(1-\rho^{2})^{q+s}
 |1-\langle\varphi(\rho \xi), \varphi(a)\rangle|^{2s}\ d\sigma(\xi) dv(u)}{(1-|a|^{2})^{-s}|1-\langle\varphi(\rho \xi),u \rangle|^{n+1+\alpha'}|1-\langle\rho \xi,a\rangle|^{2s}}
\\
 &\;&
\lesssim (1-\rho^{2})^{q+s}\int_{0}^{1}(1-r^{2})^{\alpha'-q-s}\int_{S_{n}}(1-r^{2})^{q}|g(r \eta)|^{p}(1-|\varphi_{\varphi(a)}(r \eta)|^{2})^{s}\\
 &\;&
\times \left\{\int_{S_{n}}\frac{(1-|a|^{2})^{s}|1-\langle\varphi(\rho \xi), \varphi(a)\rangle|^{2s}\ d\sigma(\xi)}{(1-
|\varphi(a)|^{2})^{s}|1-\langle\varphi(\rho\xi), r\eta\rangle|^{n+1+\alpha'}|1-\langle\rho \xi,a\rangle|^{2s}}\right\}d\sigma(\eta)dr\\
 &\;&
\lesssim (1-\rho)^{q+s}||g||_{p,q,s}^{p}\int_{0}^{1}\frac{(1-r)^{\alpha'-q-s}}{(1-r\rho)^{\alpha'+1}}\ dr
\asymp ||g||_{p,q,s}^{p}.
\end{eqnarray*}
This means that $C_{\varphi}$ is a bounded operator on $H^{p,q,s}(B_{n})$.

\vskip2mm
(6) \  Let $\{h_{k}\}$ be  any sequence  which converges to 0 uniformly on any compact set of $B_{n}$ and $||h_{k}||_{p,q,s}\leq 1$ for all $k\in \{1,2,\cdots\}$. Let $1/2<\rho<1$ and
$$
M_{\rho}=\sup_{\rho<|z|<1}\frac{1-|z|^{2}}{1-|\varphi(z)|^{2}} \ \ (z\in B_{n}).
$$

 When $p>1$,  by (3.5), H\"{o}lder inequality,  $(1-|r\xi|^{2})/(1-|\varphi(r\xi)|^{2})\leq M_{\rho}$ for any $\rho<r<1$,
  Fubini theorem,   Lemma 1.8 in [24], (3.3), Lemma 6 in [27], we have that
\begin{eqnarray*}
&\;& \sup_{0\leq r<1}(1-r^{2})^{q}\int_{S_{n}}|h_{k}[\varphi(r\xi)]|^{p}(1-|\varphi_{a}(r\xi)|^{2})^{s}\
d\sigma(\xi)\lesssim \max_{|u|\leq \rho}|h_{k}[\varphi(u)]|^{p}
\\
&\;&
+\sup_{\rho<r<1}(1-r^{2})^{q}\int_{S_{n}}(1-|\varphi_{a}(r\xi)|^{2})^{s}|G_{a,h_{k}}[\varphi(r\xi)]|^{p}|1-\langle\varphi(r\xi), \varphi(a)\rangle|^{2s}\
d\sigma(\xi)\\
&\;&
\lesssim \max_{|u|\leq \rho}|h_{k}[\varphi(u)]|^{p}+M_{\rho}^{\alpha_{1}-\alpha}\sup_{\rho<r<1}(1-r^{2})^{q+s+\alpha-\alpha_{1}}(1-|a|^{2})^{s}
\\
&\;&
\times
\int_{B_{n}}\frac{|h_{k}(u)|^{p}(1-|u|^{2})^{\alpha_{1}}}{|1-\langle\varphi(a),u\rangle|^{2s}}\left\{\int_{S_{n}}\frac{|1
-\langle\varphi(r\xi), \varphi(a)\rangle|^{2s}
\ d\sigma(\xi)}{|1-\langle\varphi(r\xi), u\rangle|^{n+1+\alpha}|1-\langle r\xi,a\rangle|^{2s}}\right\}dv(u)
\\
&\;&
\lesssim \max_{|u|\leq \rho}|h_{k}[\varphi(u)]|^{p}+M_{\rho}^{\alpha_{1}-\alpha}\sup_{\rho<r<1}(1-r^{2})^{q+s+\alpha-\alpha_{1}}
\\
&\;&
\times
\int_{0}^{1}\frac{(1-t)^{\alpha_{1}-q-s}}{(1-rt)^{\alpha+1}}\left\{\int_{S_{n}}(1-t^{2})^{q}|h_{k}(t\eta)|^{p}(1-|\varphi_{\varphi(a)}(t\eta)|^{2})^{s}\ d\sigma(\eta)\right\}dt
\\
&\;&
\lesssim \max_{|u|\leq \rho}|h_{k}[\varphi(u)]|^{p}+M_{\rho}^{\alpha_{1}-\alpha}\rightarrow M_{\rho}^{\alpha_{1}-\alpha} \ \ ( k\rightarrow\infty).
\end{eqnarray*}

When $0<p\leq 1$, we choose $0<\varepsilon<\min\{q+s, 1+\alpha'-q-\max\{s,pn\}\}$.  By (3.5), Lemma 2.15 in [24], (3.3), Lemma 1.8 in [24], Lemma 6 in [27], we have that
\begin{eqnarray*}
 &\;&
\sup_{0\leq r<1}(1-r^{2})^{q}\int_{S_{n}}|h_{k}[\varphi(r\xi)]|^{p}(1-|\varphi_{a}(r\xi)|^{2})^{s}\
d\sigma(\xi)\\
&\;&
\lesssim \sup_{0\leq r\leq \rho}(1-r^{2})^{q}\int_{S_{n}}|h_{k}[\varphi(r \xi)]|^{p}(1-|\varphi_{a}(r\xi)|^{2})^{s}\
d\sigma(\xi)+(1-|a|^{2})^{s}\\
&\;& \times\sup_{\rho< r<1}
\int_{B_{n}}\frac{|h_{k}(u)|^{p}(1-|u|^{2})^{\alpha'}}{|1-\overline{\varphi(a)}u|^{2s}}\int_{S_{n}}\frac{
(1-r^{2})^{q+s}|1-\langle\varphi(r\xi), \varphi(a)\rangle|^{2s} \ d\sigma(\xi)}{|1- \langle u, \varphi(r\xi)\rangle|^{n+1+\alpha'}|1-\langle r\xi,a\rangle|^{2s}}dv(u)
\\
&\;&\lesssim \sup_{|w|\leq \rho}|h_{k}[\varphi(w)]|^{p}+\sup_{\rho<r<1}M_{\rho}^{\varepsilon}(1-r^{2})^{q+s-\varepsilon}(1-|a|^{2})^{s}\\
&\;&
\times \
\int_{B_{n}}\frac{|h_{k}(u)|^{p}(1-|u|^{2})^{\alpha'}}{|1-\langle u, \varphi(a)\rangle|^{2s}}\left\{\int_{S_{n}}\frac{
|1-\langle\varphi(r\xi), \varphi(a)\rangle|^{2s} \ d\sigma(\xi)}{|1- \langle u, \varphi(r\xi)\rangle|^{n+1+\alpha'-\varepsilon}|1-\langle r\xi, a\rangle|^{2s}}\right\}dv(u)
\\
&\;&\lesssim \sup_{|w|\leq \rho}|h_{k}[\varphi(w)]|^{p}
+M_{\rho}^{\varepsilon}\sup_{\rho<r<1}(1-r)^{q+s-\varepsilon}||h_{k}||_{p,q,s}^{p}
\int_{0}^{1}\frac{(1-t)^{\alpha'-q-s}}{(1-rt)^{1+\alpha'-\varepsilon}} \ dt\\
&\;&\lesssim \sup_{|w|\leq \rho}|h_{k}[\varphi(w)]|^{p}+M_{\rho}^{\varepsilon}\rightarrow M_{\rho}^{\varepsilon} \ \ (k\rightarrow\infty).
\end{eqnarray*}

It follows from (3.1) that $C_{\varphi}$ is a compact operator on $H^{p,q,s}(B_{n})$.
\vskip2mm
This proof is complete. \ \ $\Box$

\vskip2mm
 {\bf Theorem 3.2} \ Let \ $p>0$, \ $0<s<n$, \ $q+s\geq 0$.  Suppose that \ $\varphi: \ B_{n}\rightarrow B_{n}$ is a holomorphic mapping.
 If $C_{\varphi}$ is a compact operator from $H^{\frac{pn}{q+n}}(B_{n})$ to $H^{p,q,s}(B_{n})$, then (3.1) holds.  If $C_{\varphi}$ is a compact operator on $H^{\frac{pn}{q+n}}(B_{n})$, then $C_{\varphi}$ is a compact operator from $H^{\frac{pn}{q+n}}(B_{n})$ to $H^{p,q,s}(B_{n})$.
 \vskip2mm
 {\bf Proof} \ If $C_{\varphi}$ is compact from $H^{\frac{pn}{q+n}}(B_{n})$ to $H^{p,q,s}(B_{n})$, then it follows from (3.4) that (3.1) holds. Conversely, $C_{\varphi}$ is a compact operator on $H^{\frac{pn}{q+n}}(B_{n})$ means that $C_{\varphi}$ is a bounded operator on $H^{\frac{pn}{q+n}}(B_{n})$. Let $\{h_{k}\}$ be  any sequence  which converges to 0 uniformly on any compact set of $B_{n}$ and $||h_{k}||_{\frac{pn}{q+n}}\leq 1$ for all $k\in \{1,2,\cdots\}$.
\vskip2mm
When  $-n<q<0$,   by  H\"{o}lder's
inequality, Lemmas 2.3, the boundedness of $C_{\varphi}$ on $H^{\frac{pn}{q+n}}(B_{n})$, we have that
\begin{align}
&
||C_{\varphi}h_{k}||_{p,q,s}^{p}\leq \sup_{0\leq r<1, a\in B_{n}}(1-r^{2})^{q+s}(1-|a|^{2})^{s} \nonumber\\
&
\times\left\{\int_{S_{n}}|h_{k}[\varphi(r\xi)]|^{\frac{pn}{q+n}}\ d\sigma(\xi)\right\}^{\frac{q+n}{n}}\left\{\int_{S_{n}}\frac{d\sigma(\xi)}{|1- \langle r\xi, a\rangle|^{\frac{2sn}{-q}}}\right\}^{\frac{-q}{n}}\nonumber\\
&
\lesssim \sup_{0\leq r<1, a\in B_{n}}\frac{(1-r^{2})^{q+s}(1-|a|^{2})^{s}||C_{\varphi}h_{k}||^{p}_{\frac{pn}{q+n}}}{(1-r^{2}|a|^{2})^{q+2s}}\lesssim ||C_{\varphi}h_{k}||^{p}_{\frac{pn}{q+n}}.
\nonumber
\end{align}

 When \ $q\geq 0$, let $p_{0}=n/(q+n)$. It follows from Lemma 2.6 that
$$
||C_{\varphi}h_{k}||_{p,q,s}^{p}
\leq
\sup_{0\leq r<1}\frac{(1-r^{2})^{q}}{(1-r^{2})^{\frac{(1-p_{0})n}{p_{0}}}}\left\{\int_{S_{n}}|h_{k}[\varphi(\sqrt{r}\xi)]|^{\frac{pn}{q+n}}\ d\sigma(\xi)\right\}^{\frac{q+n}{n}}= ||C_{\varphi}h_{k}||_{\frac{pn}{q+n}}^{p}.
$$

If $C_{\varphi}$ is  compact on $H^{\frac{pn}{q+n}}(B_{n})$,  then we have that $$\displaystyle{\lim_{k\rightarrow\infty}}||C_{\varphi}h_{k}||_{\frac{pn}{q+n}}=0 \ \Rightarrow \ \lim_{k\rightarrow\infty}||C_{\varphi}h_{k}||_{p,q,s}=0.$$
This shows that  $C_{\varphi}$ is compact from $H^{\frac{pn}{q+n}}(B_{n})$ to $H^{p,q,s}(B_{n})$.
\vskip2mm
This proof is complete. \ \ $\Box$
\vskip2mm
{\bf Statement:}  No potential conflict of interest.

\end{document}